\documentclass{article}
\usepackage{pstricks,pst-eucl}
\usepackage{pstricks-add}
\usepackage{graphicx}
\usepackage{amsthm}
\usepackage{amssymb}
\usepackage{amsmath}
\usepackage{comment}
\usepackage{amsthm}
\usepackage{color}
\usepackage{cite}
\usepackage{url}
\newtheorem{theorem}{Theorem}
\newtheorem{lemma}[theorem]{Lemma}
\newtheorem{remark}[theorem]{Remark}
\newtheorem{definition}[theorem]{Definition}
\newtheorem{proposition}[theorem]{Proposition}
\title{Isometric point-circle configurations on surfaces from uniform maps}
\author{Milagros Izquierdo$^{1}$ and Klara Stokes$^{1}$\\$^{1}$ Department of Mathematics, Link\"oping University, \\58183  Link\"oping Sweden, \\milagros.izquierdo@liu.se\\
$^{2}$ School of Engineering Sciences, University of Sk\"ovde, \\54128 Sk\"ovde Sweden, \\klara.stokes@his.se}
\date{}
\begin{document}
\maketitle

\begin{abstract}
We embed neighborhood geometries of graphs on surfaces as point-circle configurations. We give examples coming from regular maps on surfaces with maximum number of automorphisms for their genus and survey geometric realization of pentagonal geometries coming from Moore graphs. An infinite family of point-circle $v_4$ configurations on $p$-gonal surfaces with two $p$-gonal morphisms is given. The image of these configuration on the sphere under the two $p$-gonal morphisms is also described.  
\end{abstract}

\section{Introduction}
Consider the (rank two) set system of points and blocks where the points are the vertices and the blocks are the neighborhoods of the vertices of an $r$-regular graph on $v$ vertices. Such set systems are called \emph{neighborhood geometries} of graphs, and were first defined in \cite{Leemans}, within a more general context. If all neighborhoods are distinct, then the system has the following two properties: 1) each vertex appears in $r$ blocks and 2) each block contains $r$ vertices. A set system, or a geometry, with properties 1) and 2) with $v$ points and $v$ blocks is classically known as a (balanced) $v_r$ \emph{configuration}. If the intersection of each pair of blocks contains at most $d$ elements, then we say that it is of \emph{combinatorial linear dimension $d$}, since two distinct linear spaces of dimension $d$ can intersect in $d$ linearly independent points, but not in $d+1$. The blocks of a (combinatorial) geometry of linear dimension $d=1, 2, 3$ will sometimes be called lines, planes and 3-spaces, etc. A combinatorial configuration is \emph{linear} if it is of linear dimension $1$. Most configurations in the literature are linear. 

A classical example of a configuration of combinatorial linear dimension $2$ is M\"obius $8_4$ configuration, which is also a geometric configuration of 8 planes intersecting in quadruples on 8 points. The literature also contains many examples of geometric point-circle configurations. A classical example is the Miquel $(8_3,6_4)$ configuration with $6$ circles intersecting in triples on $8$ points. Two circles in the real plane intersect in at most 2 points, so combinatorially non-linear point-circle configurations in the plane correspond to configurations of combinatorial linear dimension 2. Note that if the point-circle configuration is embedded on a surface of genus $g>0$, then two circles may intersect in more than 2 points. 

The neighborhood geometry of a graph will be linear exactly when the graph does not contain any cycle of length 4. 
There is a combinatorial polarity (i.e. a duality of order two) in a neighborhood geometry defined by mapping each vertex to its neighborhood. 
If the graph is bipartite, then it defines two disconnected neighborhood geometries. The polarity then maps a point in the first connected component to its neighborhood which is then a block in the second component. Therefore, the two components are duals of each other, i.e. one is obtained from the other by interchanging the roles of the points and the blocks. This is true also if the graph is not $r$-regular. 
If the graph is not bipartite then the set system it defines consists of a single connected component, and so it is self-dual. The two disjoint geometries defined by a bipartite $r$-regular graph are dual but are not necessarily isomorphic and then not self-dual, although they have the same parameters. This is the case for example if the graph is the incidence graph of a configuration which is not self-dual. 

All the properties described above can be found in \cite{Leemans}, although it is there erroneously stated that these geometries are always self-dual (it is true for example when the graph is not bipartite and when the automorphism group of the graph is arc-transitive).

Combinatorially, the incidence graph of the neighborhood geometry of a graph is the Kronecker double cover of the graph \cite{Pisanski}. Any combinatorial property of neighborhood geometries of graphs is therefore a property of the Kronecker cover of graphs. 

The neighborhood geometry of a graph has later been applied to 1-skeletons of regular polytopes in real Euclidean $d$-space in order to construct \emph{geometric} point-hyperplane realizations (i.e. of linear dimension $d-1$) of self-polar symmetric configurations (symmetric as in \emph{with the maximal number of symmetries})\cite{Gevay}. The construction was also generalized there by using the $t$-neighborhoods of the vertices of the polytope, that is, the vertices at distance $t$ from a given vertex. 

It was observed in \cite{Gevay} that the neighborhood geometry of the 1-skeleton of a spherical polytope in real Euclidean 3-space, which is a point-plane configuration, also defines a point-circle configuration in the real Euclidean plane through stereographic projection, whenever the points in each plane are concyclic. 
Indeed, the circle-preserving property of the stereographic projection implies that any point-circle configuration drawn on the sphere can also be drawn in the real Euclidean plane. 

In \cite{Pisanski}, point-circle configurations were realized in the plane as neighborhood geometries of unit-distance graphs. 
For the neighborhood geometry of a graph embedded in the plane to be realized in terms of the circles passing through all the vertices in each neighborhood of the graph, it is necessary to ensure that these vertices are concyclic. 
Since three points define a circle, if the graph is 3-regular this condition is always satisfied. 
If the valencies of the vertices is larger than three, then embeddings with concyclic neighborhoods are special. A unit-distance embedding of the graph in the plane is an example of an embedding with this property. The 1-skeleton of a quasi-regular polyhedron in Euclidean 3-space is an example of an embedding with the same property in three dimensions.  

In \cite{StokesIzquierdo}, the construction of point-circle configurations on spherical polyhedron was generalized to surfaces in general. The motivation was there to give geometric realizations as point-circle configurations of certain pentagonal geometries coming from Moore graphs.

For a more general overview of lineal and point-circle configurations, see the books \cite{Grunbaum,PiSe}. 


In this article we will use the construction from \cite{StokesIzquierdo} to construct point-circle realizations of neighborhood geometries on several surfaces of interest. 

\section{Constructing configurations of points and isometric circles on surfaces}
Let $U$ be either the Riemann sphere, the complex Euclidean plane or the hyperbolic plane. 
A uniform tiling of $U$ is a collection of congruent polygons that partitions and fills up the entire space. 
If this tiling has $p$ $q$-gons meeting in each vertex, then the stabilizer of the tiling is a cocompact subgroup $G$ of a triangle group $\Gamma(p,2,q)$ or, if we allow orientation-reversing elements, a discrete torsion-free group of automorphisms of $U$ in which  $\Gamma(p,2,q)$ has index 2. 
Since the polygons are congruent, the neighbors of each vertex are concyclic on isometric circles. 
The distance is the spherical, the Euclidean or the hyperbolic distance respectively.  

An (compact) orientable Riemann surface is a closed topological surface $S$ with analytic structure. 
A non-orientable Riemann surface is a closed topological surface $S$ with dianalytic structure where the conjugation $z\rightarrow \bar{z}$ is allowed.
The quotient of $U$ under the action of $G$ is a Riemann surface $S=U/G$. 
The group is called the surface group of $S$ and $U$ is its universal covering space. 
By the Poincar\'e uniformization theorem any Riemann surface is the quotient $S=U/G$, where $U$ is either the Riemann sphere, the complex Euclidean plane or  the hyperbolic plane and $G$ is a discrete torsion-free group of automorphisms of $U$, possibly with orientation-reversing elements. 
The quotient of the polygonal tiling by the action of $G$ is a uniform map of type $\{p,q\}$ on the surface~\cite{Conder,Gonzalez,JS,SS}. 
In a uniform map of type $\{p,q\}$ the vertices have valency $p$, the edges have valency 2 and the faces have valency $q$. 
A map is regular if its automorphism group acts transitively on the triples of incident vertices, edges and faces, so a regular map is always uniform. 

The image of the isometric circles through the neigbhborhoods of the vertices of a uniform tiling of $U$ under the quotient by $G$ are  isometric circles through the neighbors of each vertex of the corresponding uniform map of $U/G$.  
The result is a configuration of a finite number of points and circles on the surface. 
Each circle contains $p$ points and $p$ circles goes through each point. 
We have proved the following result.    
\begin{theorem}
\label{thm1}
\cite{StokesIzquierdo} A uniform map on a surface produces a configuration of points and isometric circles on the same surface. 
\end{theorem}

Since the map completely determines the geometric point-circle configuration, the automorphism group of the configuration coincides with the automorphism group of the map. Therefore our construction will give configurations with many geometric symmetries when applied to regular maps. This motivates the study of point-circle configurations defined by regular maps in general.

Two non-isomorphic graphs can define the same neighborhood geometry. For example, both the Petersen graph and the Desargues graph has the Desargues configuration as neighborhood geometry. However, since the Desargues graph is bipartite (it is the incidence graph of the Desargues configuration), it defines the Desargues configuration as a point-circle configuration twice. 
Also, the same $r$-regular graph can be embedded in a Riemann surface as a uniform map in several ways. Consequently, there may be many ways to realize the same configuration in terms of points and circles on some surface.  

In a paper from 1949 Coxeter explored the relation between self-dual configurations and arc-transitive graphs (he used the term \emph{regular graph}) \cite{Coxeter}. He embedded the incidence graph of the configurations as regular maps on surfaces. 
For bipartite graphs, this is exactly what we also do. However,  our approach goes further. We obtain a geometric configuration defined by the incidences of elements from two classes of distinct geometric objects, points and circles, on the surface. As an immediate consequence of this, we see that all configurations represented in \cite{Coxeter} as regular maps of incidence graphs, are actually point-circle configurations on the same surfaces.  

For graphs that are not bipartite, our construction is essentially different from Coxeter's approach. However, his general principle \emph{``interesting configurations are represented by interesting graphs''} \cite{Coxeter} may still be applied.   
\section{Point-circle configurations from some classical regular maps}

In this section we apply Theorem \ref{thm1} to some classic regular maps and give the resulting point-circle configuration. 

\subsection{Klein's quartic}
Klein's quartic projective curve, given by the equation $x^3y+y^3z+z^3x=0$ over the complex field is the curve of smallest genus that attains the Hurwitz bound. Its genus is $g=3$ and its automorphism group is $PSL(2,7)$ of order $84(g-1)=168$, so it is the curve of genus 3 with maximum number of automorphisms.  There is an epimorphism from the triangle group $\Gamma(2,3,7)$ to $PSL(2,7)$ and its kernel is the surface group uniformizing Klein's quartic, as Riemann surface. This results in the classical regular heptagonal map on the surface of type $\{3,7\}$ with 56 vertices of valency 3, 84 edges and 24 heptagonal faces. The neighborhood geometry of the underlying non-bipartite graph is an autopolar $56_3$ configuration, which is linear, since the graph has no 4-cycles. 

The dual map (obtained by interchanging the roles of vertices and faces) has 24 vertices of valency 7, 84 edges and 56 triangular faces.  The 24 vertices corresponds to the 24 Weierstrass points of the surface. 
The neighborhood geometry of the underlying non-bipartite graph is an autopolar $24_7$ configuration of combinatorial linear dimension 2. Given any point $p$, there are exactly two points which are not concylic with $p$. 
All other points are concyclic with $p$ exactly twice. 








\subsection{Bring's curve}
Bring's curve is a complex projective curve of genus 4 given by the equations $\sum_{i=1}^5x_i=\sum_{i=1}^5x_i^2=\sum_{i=1}^5x_i^3=0$. Its automorphism group is the symmetric group acting on 5 elements. It is the curve of genus 4 with maximum number of automorphisms.

Consider the triangle group $\Gamma(2,4,5)$. There is an epimorphism from $\Gamma(2,4,5)$ to the symmetric group $S_5$ and the kernel is the surface group of Bring's curve (as Riemann surface), which is normal in $\Gamma(2,4,5)$. The surface allows a regular map of type $\{4,5\}$. This map has 30 vertices of valency 4, 60 edges and 24 pentagonal faces. The neighborhood geometry of the underlying graph is an autopolar $30_4$ configuration. Since the graph has no 4-cycles, the configuration is linear.   

The dual map has 24 vertices, 60 edges and 30 quadratic faces. The underlying graph is bipartite, and is therefore the incidence graph of its neighborhood geometry, an autopolar $12_5$ configuration of combinatorial linear dimension 2.  
 Given any point $p$, there is exactly one point which is not concylic with $p$. 
All other points are concyclic with $p$ exactly twice. 
This gives the configuration an antipodal property.





\subsection{The Bolza curve}
The Bolza curve is a complex projective curve of genus 2 with automorphism group $GL(2,3)$ of order 48. It is the curve of genus 2 with maximum number of automorphisms. Its surface group (as Riemann surface) is a normal subgroup of the triangle group $\Gamma(2,3,8)$. It has a regular map of type $\{3,8\}$ with 16 vertices, 24 edges and 6 octagonal faces. 
The neighborhood geometry of the underlying bipartite graph is the unique linear $8_3$ configuration known as the M\"obius-Kantor configuration, which is therefore realized as a point-circle configuration on the Bolza curve. This configuration has the antipodal property, mutually non-collinear points occur in pairs.  
The dual map has 6 vertices, 24 edges and 16 triangular faces. The neighborhood geometry of the underlying graph is a degenerated $6_4$ configuration of 3 circles through 6 points where each circle appears twice. Each pair of distinct circles meet in 2 points.




  
\section{Pentagonal geometries as point-circle configurations from Moore graphs}
A pentagonal geometry is a (linear) combinatorial configuration with the property that, for any point $p$, all points that are not collinear with $p$ are on a single line, which is called the opposite line of $p$~\cite{Ball}. 
The lines in a pentagonal geometry are of two types, lines that are the opposite line of some point, and lines that are not.  

A pentagonal geometry in which all lines are opposite lines is self-polar by the polarity that make correspond each point to its opposite line. 
The reduced Levi graph  (in the sense of \cite{Artzy}) defined by the polarity of a self-polar pentagonal geometry is the graph in which the vertices are pairs of one point and its polar line and two vertices are joined by an edge if the point of one vertex is incident with the line of the other vertex. Therefore this graph is exactly the deficiency graph of the geometry, that is, the graph in which the vertices are the points and two vertices are joined by an edge if the points are not collinear. The pentagonal geometry can be recovered from the reduced Levi graph as its neighborhood geometry. 
More generally, the neighborhood geometry of the reduced Levi graph of an autopolar configuration is always equal to the original configuration. 
This construction of pentagonal geometries was first described in \cite{Ball}, where it also was proved that pentagonal geometries with $r=k$ are exactly the ones with a Moore graph of diameter two as reduced Levi graph. 

There are only three known Moore graphs of diameter 2, the cycle graph of length 5, the Petersen graph and the Hoffman-Singleton graph. 
These graphs have degree 2, 3 and 7, respectively. 
The existence of a Moore graph of degree 57 is still an open question. 
The pentagonal geometries obtained from these graphs are, respectively, the ordinary pentagon, the Desargues' configuration and a pentagonal geometry with parameters $(7,7)$ and with 50 points and 50 lines. 
In \cite{Ball}, it was also proved that all pentagonal configurations of order $(k,k+1)$ can be constructed from pentagonal geometries of order $(k+1,k+1)$ through the removal of one point and its opposite line. 
There are therefore at most three such pentagonal geometries, with $k=2,6$ and maybe $56$. 

Regular embeddings of the two smallest Moore graphs are well-known. The cycle graph can be embedded with full automorphism group as a regular map on the Riemann sphere with two pentagonal faces. The Petersen graph has an embedding in the real projective plane as a regular map with six pentagonal faces, obtained from the classical spherical map $\{3,5\}$ by identifying antipodal points. By Theorem \ref{thm1}, this implies that both the pentagon and the Desargues configuration allows realizations in terms of points and isometric circles with full automorphism group on surfaces of orientable genus 0 and non-orientable genus 1, respectively. 

There is no regular embedding of the Hoffman-Singleton graph, but there are uniform pentagonal embeddings of type $\{7,5\}$ on non-orientable surfaces of genus 57 with automorphism group of the map either trivial, of order 5 or of order 7 \cite{ConderStokes}. By Theorem \ref{thm1} this implies that the pentagonal geometry $(7,7)$ on 50 points and 50 blocks can be realized as a point-circle configuration on a surface of non-orientable genus 57 with any of these three automorphism groups \cite{StokesIzquierdo}.  It is also possible to realize this pentagonal geometry as a point-hypersphere configuration with full automorphism group in Euclidean space of 24 dimensions as the geometric neighborhood geometry of the well-known embedding of the Hoffman-Singleton graph in the Leech lattice \cite{StokesIzquierdo}.

\section{Point-circle configurations on $p$-gonal surfaces with two cyclic $p$-gonal morphisms}
The universal covering space $U$ of an (orientable) Riemann surface $U/G$ covers it with infinitely many sheets. Each point of $U/G$ is the representative of exactly one orbit (fiber) of the points in $U$ under the action of the surface group $G$, which is a torsion-free discrete group of automorphisms of $U$, possibly with orientation-reversing elements. 

By considering also cocompact discrete groups $H$ of automorphisms of $U$ with elliptic elements, one obtains a surface $U/H$ with singular points, a geometric orbifold. An orbifold is a more general concept than a Riemann surface. When $H$ has no elliptic elements, then $U/H$ is a Riemann surface. 
 
Let $G$ be a finite index $n$ subgroup of a discrete subgroup $H$ of automorphisms of $U$ (a Fuchsian group). The inclusion $G\hookrightarrow H$ induces a (possibly ramified) covering $f:U/G\rightarrow U/H$ of degree $n$. The covering $f$ is determined by the action of $H$ on the $G$-cosets $\theta:H\rightarrow \Sigma_{|H:G|}$. If $G$ is normal in $H$, then the covering $f:U/G\rightarrow U/H$ is a regular covering given by the monodromy $\theta:H\rightarrow H/G$. 
Assume that $G$ has no elliptic elements, so that $G$ is a surface group uniformizing a Riemann surface $U/G$. Assume also that the covering morphism $U/G\rightarrow U/H$ of degree a prime number $p$. Then there is an automorphism $\phi:S\rightarrow S$ such that the deck-transformation group of the covering is generated by $\phi$, that is, $\langle \phi \rangle \sim H/G=C_p$. 
If the genus of the underlying surface of $U/H$ is 0, then we say that the surface $U/G$ is a $p$-gonal surface and $f:U/G \rightarrow U/H$ is a $p$-gonal morphism \cite{CoIzYi,Gonzalez2}. 
 

According to Castelnuovo-Severi \cite{Accola} a compact Riemann surface which allows more than one $p$-gonal morphism has genus $g$ satisfying $g\leq (p-1)^2$. 
Additionally, it is known that if $S_g$ has several $p$-gonal morphisms, then these morphisms are all conjugate \cite{Gonzalez2}. 

For all prime $p\geq 3$ there is a family of surfaces with two distinct cyclic $p$-gonal morphisms \cite{CoIzYi} of genus $(p-1)^2$, implying that the Castelnuovo-Severi inequality is sharp. 
For each $p$, one of the surfaces in this family, which we call $Y_p$, is quasi-platonic, i.e. its surface group is normal in the triangle group $\Gamma(4,2,2p)$. 
The automorphism group of this surface\footnote{In \cite{CoIzYitrigonal} an erroneous presentation of this group for $p=3$ was given. A presentation is  $ \langle a, b, s, t / a^3=b^3=s^2=t^4 =(st)^2 =(sa)^2=sbsb^{2}= t^3atb^2=t^3bta=1 \rangle$, as correctly stated in \cite{Ying}. The two $p$-gonal automorphisms are then $ab$ and $ab^{-1}$} is $(C_p\times C_p)\rtimes D_4$, where $C_n$ and $D_n$ are the cyclic and the dihedral groups of order $n$ and $2n$, respectively.   
Since the surface group is a normal subgroup in $\Gamma(4,2,2p)$, the surface allows a regular map of type $\{4,2p\}$ of genus $(p-1)^2$. 
The underlying graph of this map is a bipartite $4$-regular symmetric graph on $2p^2$ vertices of girth 4. 
The neighborhood geometry of this graph has $p^2$ points and $p^2$ blocks and it is self-polar, because the graph is symmetric. Each pair of points appears in exactly 0, 1 or 2 blocks, so it is not linear, but (combinatorially) planar. 
By Theorem \ref{thm1}, the geometry is realized as a point-circle configuration on the surface $Y_p$. There are 4 points on each circle and 4 circles through each point.

The neighborhood geometry of the underlying bipartite graph of the dual map is less interesting, since it is a degenerate $2p_{2p}$ configuration of one circle through $2p$ points, for the $2p$ circles are all the same.  

The automorphism group of $Y_p$ contains two conjugate $p$-gonal morphisms such that the orbit space of their action on $Y_p$ is the Riemann sphere with $2p$ singular points, all of degree $p$. 
The action of each one of the two $p$-gonal morphism on the regular map of type $\{4,2p\}$ on $Y_p$ divides the vertices in $2p$ orbits of $p$ vertices each. The result is a $2$-regular bipartite graph on $2p$ vertices embedded as a map on the sphere. 
By Theorem \ref{thm1}, its neighborhood geometry is a point-circle configuration with $p$ points and $p$ circles on the sphere,  2 points on each circle and 2 circles through each point. 
We have proved the following. 
\begin{theorem}
There is an infinite family of combinatorially planar auto-polar configurations $p^2_4$, realizable as point-circle configurations on orientable $p$-gonal surfaces of genus $(p-1)^2$. The orbits of such a configuration under the action of the two $p$-gonal morphisms defines a $p_2$ configuration of $p$ circles and $p$ points on the sphere which is isomorphic to the cycle graph of order $p$. 
\end{theorem}

The smallest member of this family of point-circle configurations on $p$-gonal surfaces is perhaps the most interesting one. 
When $p=3$, then the graph of the regular map is bipartite on 18 vertices. Its neighborhood geometry has 9 points and 9 blocks.  Each pair of points are in at least one block, therefore either in 1 or 2 blocks. Interestingly, given a point $x$, the points that are not in the same block twice (and therefore once) with $x$ form a single block. This property can be seen as a generalization of the property defining a pentagonal (linear) geometry: given a point $x$, the points that are not collineal with $x$ are collinear on a single line of the geometry, and they are the only points on that line. 
The property in its general form for a configuration of combinatorial linear dimension $d$ is as follows. 
\begin{definition}
A configuration of combinatorial linear dimension $d$ has the generalized pentagonal property if, given any point $x$, there are points which are in the same block as $x$ exactly $d$ or $d-1$ times and the points that are in the same block as $x$ exactly $d-1$ times form a single block. 
\end{definition}
For $d=1$ this is the defining property of the pentagonal geometries. For $d=2$, there is the $9_4$ configuration we just described. Its incidence graph, the 4-regular symmetric graph on 18 vertices is listed as the second symmetric graph of order 18 in \cite{Marstonlist}. In the same list one also finds an example for $d=3$, the neighborhood geometry defined by the third symmetric graph of order 26, a 6-regular symmetric graph on 26 vertices with automorphism group of order 156. These examples are all auto-polar configurations by construction and somehow they generalize the pentagon. A $d$-dimensional auto-polar geometry with the generalized pentagonal property must have  $1+r^2/d$ points and $1+r^2/d$ blocks. In particular $d$ must divide the number of points per block $r$. The pentagonal geometries are in general not auto-polar linear configurations and may be non-balanced, that is, they may have more lines through a point than points on a line. It would be interesting to know if there are examples of non-balanced configurations of combinatorial linear dimension $d$ with the generalized pentagonal property. 

The neighborhood geometry of the underlying symmetric graph of the regular $\{4,6\}$ map on $Y_3$ is also the smallest member of the family of planar point-circle configurations defined as the neighborhood geometries of the generalized cuboctahedron graph in \cite{Pisanski}. A part from this example, the family of point-circle configurations from generalized cuboctahedron graphs and the family defined by the regular maps of type $\{4,2p\}$ on the surfaces $Y_p$ are disjoint. In \cite{Pisanski} it is stated as an open question whether there is an isometric realization in terms of points and circles of any of the configurations in the family coming from the generalized cuboctahedron graphs. In this article we have given a realization of the first member of that family in terms of points and isometric circles on an orientable surface of genus 4. 


This special configuration can also be realized as a point-circle configuration in the complex plane as the neighborhood geometry of the $3\{4\}2$ regular and complex polygon.  The points and edges of this polygon are the points and lines of the generalized quadrangle with parameters $(2,1)$ (a grid on 9 points with 3 points on each line and 2 lines through each point). In \cite{Leemans} neighborhood geometries were defined for geometries as well as for graphs. Note that the neighborhood geometry of the generalized quadrangle (as a geometry) equals the neighborhood geometry of its symmetric, distance-regular collinearity graph obtained by forgetting the geometric structure. 

It is described in \cite{Coxeterbook} how a complex polygon can be represented in the real plane, and a real representation of this particular $3\{4\}2$ polygon can be found on page 108. The regularity of the polygon and the properties of this representation imply that the graph that has as vertices the vertices of the real representation and as edges the sides of the triangles representing the complex edges of the polygon in the real plane, is a unit-distance graph in the real plane.  
This implies that this particular neighborhood geometry can be realized as an isometric point-circle configuration in the Euclidean real plane. This gives also a planar answer to the open question mentioned before from \cite{Pisanski}. 

\section*{Acknowledgements}
The authors want to thank Marston Conder for helpful discussions. Calculations were mainly done with Magma. The second author acknowledges partial financial support from the Spanish MEC project ICWT (TIN2012-32757).

\end{document}